\documentclass{amsproc}

\usepackage{amsmath}
\usepackage{amssymb}
\usepackage{amscd}
\usepackage{amsthm}

\def\Zee{\mathbb{Z}}
\def\Q{\mathbb{Q}}
\def\Ar{\mathbb{R}}

\def\Ker{\operatorname{Ker}}

\def\Hom{\operatorname{Hom}}

\def\Id{\operatorname{Id}}

\newtheorem{theorem}{Theorem}[section]
\newtheorem{lemma}[theorem]{Lemma}
\newtheorem{corollary}[theorem]{Corollary}
\newtheorem{prop}[theorem]{Proposition}
\newtheorem{conj}[theorem]{Conjecture}
\theoremstyle{definition}

\newtheorem{example}[theorem]{Example}

\theoremstyle{remark}
\newtheorem{remark}[theorem]{Remark}

\numberwithin{equation}{section}

\begin{document}

\title[Mirror
Symmetry and Integral Variations of Hodge Structure]{Mirror
Symmetry and Integral Variations of Hodge Structure Underlying One
Parameter Families of Calabi-Yau Threefolds}

\author{Charles F. Doran}
\address{Department of Mathematics, University of Washington, Seattle,
Washington 98195}
%\curraddr{Department of Mathematics and Statistics, Case Western
%Reserve University, Cleveland, Ohio 43403}
\email{doran@math.washington.edu}

\author{John W. Morgan}
\address{Department of Mathematics, Columbia University, New York, New York
10027}
\email{jm@math.columbia.edu}

%    General info
%\date{January 1, 1994 and, in revised form, June 22, 1994.}


%\keywords{Differential geometry, algebraic geometry}

\begin{abstract}
This proceedings note introduces aspects of the authors' work
relating mirror symmetry and integral variations of Hodge
structure.  The emphasis is on their classification of the
integral variations of Hodge structure which can underly families
of Calabi-Yau threefolds over ${\mathbb{P}}^1 \setminus
\{0,1,\infty\}$ with $b^3 = 4$, or equivalently $h^{2,1} = 1$, and
the related issues of geometric realization of these variations.
The presentation parallels that of the first author's talk at the
BIRS workshop.
\end{abstract}

\maketitle

\vspace{.1in}

\section{Integral structures in mirror symmetry}

\subsection{The first examples}\label{firstegs}

Since its introduction to the mathematical community through the
seminal papers of Greene-Plesser \cite{GP} and Candelas-de la
Ossa-Green-Parkes \cite{CdOGP}, Mirror Symmetry has been the source of
the most persistently rich and subtle novel mathematics yet to
emerge from the study of string dualities.  One of its key
features is its resistance to rigorous mathematical definition and
even more to being described within any single traditional
mathematical setting. Research in mathematics related to mirror
symmetry is thus driven by the goal of discerning deeper, more
complete, and purely mathematical avatars of this very physical
duality.

For example, one of the earliest mirror symmetric mathematical
predictions was that a ``mirror pair'' of Calabi-Yau threefolds
$X$ and $\widetilde{X}$ should have the property that their Hodge
numbers ``mirror'' one another, i.e., that  $h^{1,1}(X) =
h^{2,1}(\widetilde{X})$ and $h^{2,1}(X) = h^{1,1}(\widetilde{X})$.

The first proposal for such a pair was introduced in \cite{GP}
where it was applied to Calabi-Yau threefolds presented as
hypersurfaces in a weighted projective space.  The construction
involves the operation of quotienting by a finite group (or
``orbifolding'').

The simplest of Calabi-Yau threefold hypersurfaces, the generic
quintic hypersurface in ${\mathbb{P}}^4$, often denoted by
${\mathbb{P}}^4[5]$, has $h^{1,1} = 1$ (coming from the
polarization class) and $h^{2,1} = 101$ (corresponding to 101
complex structure moduli). This last follows from Kodaira-Spencer
theory \cite{Kod} and the theorem of Bogomolov-Tian-Todorov
\cite{Tian} that any Calabi-Yau manifold admits a locally
universal deformation over a smooth base. In this case all of the
complex structure deformations arise from varying the coefficients
of the defining equations for quintic hypersurfaces in
${\mathbb{P}}^4$. Let $G$ denote the finite group
$({\mathbb{Z}}/5{\mathbb{Z}})^3$ presented as
$$\left\{(a_1,a_2,a_3,a_4,a_5) \in ({\mathbb{Z}}/5{\mathbb{Z}})^5
\ \left\bracevert \ \sum_{i=1}^5 a_i = 0 \right. \right\} /
({\mathbb{Z}}/5{\mathbb{Z}}) \ ,$$ where the action of
${\mathbb{Z}}/5{\mathbb{Z}}$ is the diagonal one. Upon quotienting
the projective space ${\mathbb{P}}^4$ by $G$, acting as
$$(x_1:x_2:x_3:x_4:x_5) \mapsto
(\mu^{a_1}x_1:\mu^{a_2}x_2:\mu^{a_3}x_3:\mu^{a_4}x_4:\mu^{a_5}x_5)
$$ (where $\mu = \exp{2 \pi \imath / 5}$), resolving the resulting
orbifold singularities, and keeping track of the hypersurface
throughout the process, another family of Calabi-Yau threefolds is
constructed --- no longer as hypersurfaces in projective space,
but now with the property that generically $h^{1,1} = 101$ and
$h^{2,1} = 1$.  It follows that these threefolds sit in a
one-parameter family of complex structures.  The specific
presentation of the mirror family $\widetilde{X}_z$ in
${\mathbb{P}}^4 / G$ is given by
$$x_1^5 + x_2^5 + x_3^5 + x_4^5 + x_5^5 - 5 z^{-1/5} x_1 x_2 x_3
x_4 x_5 = 0 \ .$$ Here the complex deformation parameter is
identified with a coordinate on the base of the family
${\mathbb{P}}^1 \setminus \{0,1,\infty\}$.

\begin{remark}
In \cite[Table 1]{GP} this construction is extended to include
intermediate quotients by subgroups of $G$.  In particular,
instead of the quintic itself one can consider the quotient of the
quintic by the free ${\mathbb{Z}}/5{\mathbb{Z}}$-action given by
$$ (x_1:x_2:x_3:x_4:x_5) \mapsto (x_1:\mu x_2: \mu^2 x_3: \mu^3
x_4: \mu^4 x_5) \ ,$$ resulting in a family of Calabi-Yau
threefolds with $h^{1,1} = 1$ and $h^{2,1} = 21$. These ``quintic
twin'' Calabi-Yau threefolds have fundamental group
${\mathbb{Z}}/5{\mathbb{Z}}$. The mirror family is constructed in
\cite{GP} as another intermediate quotient, and consists (as
expected) of  simply connected Calabi-Yau threefolds with $h^{1,1}
= 21$ and $h^{2,1} = 1$. Moreover the parameter space of this
quintic twin mirror family $\widetilde{X}^{twin}_z$ is again
${\mathbb{P}}^1 \setminus \{0,1,\infty\}$ as for the original
quintic mirror family.
\end{remark}

Since the mirror quintic family $\widetilde{\mathfrak X} = \{
\widetilde{X}_z \}$ depends on just one parameter, an old method
of Griffiths-Dwork \cite[\S 5.3]{CoKa} can be readily applied to
determine the ordinary differential equation satisfied by the
periods of the (unique up to complex scaling) holomorphic
three-form on the Calabi-Yau threefold. More formally, consider
the local system denoted $R^3 \varphi_* {\mathbb{C}}$ arising from
the third cohomology (with complex coefficients) of the fibers
$\widetilde{X}_z$ of the algebraic family $\varphi \colon
\widetilde{\mathfrak X} \rightarrow B := {\mathbb{P}}^1 \setminus
\{0,1,\infty\}$. Integration over $3$-cycles defines a pairing,
making this local system dual to the local system
 whose fiber over $z$ is $H_3(\widetilde{X}_z, {\mathbb{C}})$.
Every class $\gamma \in H_3(\widetilde{X}_z, {\mathbb{C}})$
extends uniquely as a multi-valued flat section of this later
local system. Consider the holomorphic bundle $R^3 \varphi_*
{\mathbb{C}} \otimes {\mathcal{O}}_B$; the local system
$R^3\varphi_*{\mathbb C}$ determines a flat structure on this
vector bundle and hence a flat connection, $\nabla_z$. The
sub-line bundle $\mathcal{H}^{3,0}$ of classes of relative
holomorphic $3$-forms is of special interest. A holomorphic
section $\varpi(z)$ of $\mathcal{H}^{3,0}$ amounts to a
holomorphic family of holomorphic $3$-forms. Generically
$\varpi(z)$ and its first three covariant derivatives,
$\nabla_z^i(\varpi)$ $(1\le i\le 3)$ are linearly independent, and
$\nabla_z^4(\varpi)$ is expressible as a linear combination of
these with coefficients meromorphic functions on the base $B$.
This relationship is the differential equation associated to the
family of cohomology classes $[\varpi(z)]$ --- the {\em
Picard-Fuchs equation}. There is also the {\em Picard-Fuchs ODE}
obtained by replacing $\Theta$ by $z \frac{d}{dz}$. An integral
lattice in the space of all solutions of this ODE is given by the
periods $\int_\gamma \varpi(z)$ as $\gamma$ varies over all
(integral) three cycles. There is  a holomorphic map of the dual
of the vector bundle of solutions of the Picard-Fuchs ODE to $R^3
\varphi_* {\mathbb{C}} \otimes {\mathcal O}_B$. This map is an
isomorphism near any regular singular point of the ODE.

For the quintic mirror family $\widetilde{\mathfrak X}$, letting
$\Theta = z \nabla_z$, for an appropriate holomorphic family
$\varpi(z)$ of holomorphic $3$-forms the Picard-Fuchs equation is
the generalized hypergeometric differential equation
$$\left[ \Theta^4 - z \left( \Theta + \frac{1}{5} \right)
\left( \Theta + \frac{2}{5} \right) \left( \Theta + \frac{3}{5}
\right) \left( \Theta + \frac{4}{5} \right) \right] \varpi(z) = 0
\ .$$ The Picard-Fuchs ODE has 3 singular points $0, 1, \infty$,
all of which are regular singular points. Let $\gamma_0$ be an
indivisible integral cycle invariant under monodromy about $z=0$.
A holomorphic solution to the Picard-Fuchs ODE given is the
hypergeometric-type series
\begin{equation*}
\int_{\gamma_0} \varpi(z) = \sum_{n \geq 0} \frac{(5n)!}{(n!)^5}
\left( \frac{z}{5^5} \right)^n \end{equation*}
 defined on the
punctured unit disk. This function extends across the origin,
reflecting the fact that it is integration of a family of
holomorphic $3$-forms over a flat family of $3$-cycles invariant
under monodromy around $z=0$. This power series is absolutely
convergent inside the unit disk.  Of course, it also admits a
multi-valued analytic continuation to the whole thrice-punctured
sphere given explicitly in terms of Meijer functions (Mellin-Barnes
type integrals) on ${\mathbb{P}}^1 \setminus \{0, 1, \infty \}$
\cite{HTF,GL}.

\begin{remark}
The story for the quintic twin mirror family is similar, with
everything working just as above.  In fact, for appropriate
choices of the family of holomorphic three-forms, the Picard-Fuchs
ODE defined as above and satisfied by the periods of the two
families of holomorphic forms on $\widetilde{X}_z$ and
$\widetilde{X}^{twin}_z$ are exactly the same. Thus, the
difference between the two families is not reflected in the
Picard-Fuchs differential equations.  It turns out that it is
reflected in the lattices spanned by the periods of the
holomorphic three-form over a basis of integral cycles. Of course,
the difference is also seen from the fact that the ranks of the
even dimensional cohomologies of the two families are different.
\end{remark}

Mirror symmetry predicts that the periods over a suitable integral
basis of three-cycles on the mirror quintic contain extremely subtle
information about the geometry of the original quintic hypersurface
in ${\mathbb{P}}^4$. In particular, mirror symmetry predicts the
number of rational curves of a given degree lying on the
hypersurface \cite{CdOGP,CoKa}. These curve counts can be read off
of an appropriate generating function built from the ``mirror map''
$q$-series --- the single-valued local inverse to the projectivized
period mapping about a regular singular point of maximal unipotent
monodromy (or, in terms of the geometry of the family of Calabi-Yau
threefolds, of maximal degeneracy). In many cases these curve counts, and
their descriptions in terms of hypergeometric functions, have
now been rigorously established \cite{LLY1,LLY2,Giv}.
We will not be pursuing this approach here; our focus is the
integral lattice structure and  monodromy rather than enumerative
geometry.

\subsection{Classical Hodge theory}

In order to study more general families of Calabi-Yau threefolds
than just the mirror quintic or mirror quintic twin hypersurfaces
we introduce here Hodge structures and their variations.  Consider
a Calabi-Yau threefold $\widetilde X$ with $h^{2,1}=1$. The
cohomology $H^3(\widetilde{X},{\mathbb{R}})$ is polarized by the
intersection form, denoted $\langle \cdot , \cdot \rangle$,  which
is unimodular and skew, and has a Hodge decomposition
$$H^3(\widetilde{X}, {\mathbb R}) \otimes_{\mathbb R} {\mathbb C} =
H^3(\widetilde{X}, {\mathbb{C}}) = H^{3,0} \oplus H^{2,1} \oplus
H^{1,2} \oplus H^{0,3}$$ with associated Hodge filtration
$$F^3 = H^{3,0} \ , \ F^2 = H^{3,0} \oplus H^{2,1} \ , \ F^1 =
H^{3,0} \oplus H^{2,1} \oplus H^{1,2} \ .$$
Furthermore, for all $i$, $0 \leq i \leq 3$ these filtrants satisfy the conditions
\begin{equation}
(F^i)^\perp = F^{4-i} \label{hodge1}
\end{equation}
and
\begin{equation}
F^i \oplus \overline{F^{4-i}} = H^3(\widetilde{X}, {\mathbb{C}})\ .
\label{hodge2}
\end{equation}
Notice we recover the Hodge decomposition from the filtration
since $F^i\cap \overline{F^{3-i}}=H^{i,3-i}$.

 For us a Calabi-Yau
threefold has holonomy all of $SU(3)$ and hence $H^1(\widetilde{X},
{\mathbb{C}}) = 0$, implying that all of $H^3$ is primitive. Thus,
the hermitian form
$$h(\omega_1, \omega_2) = \imath \int_M \omega_1 \wedge
\overline{\omega}_2 = \imath \langle \omega_1 ,
\overline{\omega}_2 \rangle$$ is positive definite on $H^{3,0}$
and negative definite on $H^{2,1}$, i.e., for $\omega \neq 0$,
\begin{equation}\label{hodge3}
h(\omega, \omega) \ \ \mbox{is} \ \begin{cases}> 0 & \mbox{when} \
\omega \in
H^{3,0} \\
     < 0 & \mbox{when} \  \omega \in H^{2,1} \, .\end{cases}
\end{equation}

Fix a rank 4 real vectorspace $V$ with a nondegenerate skew form.
Weight-three Hodge structures with $h^{2,1} = 1$ on $V$
(and by this we shall always mean primitive and polarized weight-three
Hodge structures)
are classified by the period domain
$${\mathcal{D}} = \mbox{Sp}(4, {\mathbb{R}}) / \mbox{U}(1)
\times \mbox{U}(1) \ .$$  ${\mathcal D}$ is embedded as a domain
inside $\mbox{Sp}(4, {\mathbb{C}})/P$, where $P$ is the maximal
parabolic subgroup. A holomorphic family of such structures on $V$
parameterized by a complex curve $U$ consists of a holomorphic
filtration by subbundles
$${\mathcal F}^3\subset
{\mathcal F}^2\subset {\mathcal F}^1\subset {\mathcal F}^0 =
V\otimes_{\mathbb R} {\mathcal O}_U,$$
 with Conditions~(\ref{hodge1}), (\ref{hodge2}), and~(\ref{hodge3}) holding on every fiber,
resulting in a family of (polarized) Hodge structures. Given such a
holomorphic family, the {\em period map} is the (holomorphic)
classifying map $\Pi \colon U \rightarrow {\mathcal{D}}$.

There is an additional geometric condition on the period map $\Pi
\colon U \to {\mathcal D}$ if the family of Hodge structures is
geometric, i.e., arises from a family of smooth varieties. This is
the {\em Griffiths transversality condition}, also called
horizontality, \cite[Condition (2.8) on page 143]{Grif} which says
that $({\mathcal F}^i)'\subset {\mathcal F}^{i-1}$ (where the
derivative is taken with respect to a local coordinate on $U$).
Horizontal families of Hodge structures are called {\em variations
of Hodge structure (VHS)}, or ${\mathbb Z}$-VHS, to emphasize  the
underlying flat integral structure. Horizontality can be
reformulated in terms of a (non-integrable) distribution, called the
{\em horizontal distribution}, on ${\mathcal D}$. Any period map
$\Pi \colon U \to {\mathcal D}$ from a horizontal family is tangent
to this distribution.

For the Hodge structures we are considering, there is the natural
map
$$p_3\colon {\mathcal D} \to {\mathbb P}(V_{\mathbb C})$$ associating to
each
Hodge filtration the line $F^3$ in $V_{\mathbb C}=V \otimes_{\mathbb
R}{\mathbb C}$. The horizontality condition for a curve $C\subset
{\mathcal D}$ is equivalent to
$$T(p_3(C))_p\subset p^\perp/\langle p\rangle\subset
(V_{\mathbb C})/\langle p\rangle=T\left({\mathbb P}(V_{\mathbb
C})\right)_p.$$ Conversely, given an analytic curve $C\subset
{\mathbb P}(V_{\mathbb C})$ satisfying this condition, there is a
unique lifting ${\mathcal C}\subset {\mathcal D}$ of $C$ that is
horizontal: Given $p\in C$, the Hodge filtration is given as follows
$F^3=\langle p\rangle$, $F^2$ is the preimage in $V_{\mathbb C}$ of
$T( C)_{ p}\subset (V_{\mathbb C})/\langle  p\rangle$, and
$F^1=(F^3)^\perp$. For a fuller explanation of this see \cite{BG}.

Now let $V=H^3(\widetilde X,{\mathbb R})$ and let $U$ be a local
simply connected open subset in the (one-dimensional) moduli space
of complex structures on $\widetilde X$. Over $U$, the family of
holomorphic $3$-forms $\varpi(z)$ is a holomorphic section of
$V_{\mathbb{C}}$. It satisfies a differential equation of the form
$$\nabla^4_z(\varpi(z)) + \sum_{i=0}^3 P_{4-i}(z) \nabla^i(\varpi(z)) = 0 \ .$$
Using a local $C^\infty$-trivialization of the family, we extend any
$3$-cycle $\gamma \in H_3(\widetilde X_{z_0},{\mathbb{C}})$  to a
family of $3$-cycles in all  neighboring fibers. The periods
$\int_\gamma \varpi(z)$ form a (local) holomorphic function on the
base. Since the family of cycles $\gamma$ is a parallel section, the
function
$$\varphi_\gamma(z) = \int_\gamma \varpi(z)$$ satisfies the
Picard Fuchs ODE (i.e., the ODE associated to the original
Picard-Fuchs equation for $\varpi$):
$$\left( \frac{d}{dz} \right)^4 \varphi_\gamma + \sum_{i=1}^4 P_{4-i}(z)
\left( \frac{d}{dz} \right)^i \varphi_\gamma = 0. $$ The
$4$-dimensional space of periods as $\gamma$ ranges over all of
$H_3(\widetilde X_{b};\mathbb{C})$ is the space of all solutions of
the Picard-Fuchs ODE. This gives an identification (over a
contractible neighborhood of $b \in B$) between the space of
solutions to the Picard-Fuchs ODE and the vector bundle $(R^3
\varphi_* {\mathbb{C}})^* \otimes {\mathcal O}_B$. Fixing a basis
$\gamma_1^*, \ldots, \gamma_4^*$ of the lattice $H^3(\widetilde
X,{\mathbb Z})/\{\rm Torsion\}$, we have
$$\varpi(z) = \sum_{i=1}^4 f^i(z) \cdot \gamma_i^* \ .$$
Of course the $f^i(z)$ are simply the
periods
$$\int_{\gamma_i} \varpi(z)$$
where $\gamma_1, \ldots, \gamma_4$ is the basis of
$H_3(\widetilde{X}_z, {\mathbb{Z}})/\{\rm Torsion\}$ dual to
$\gamma_i^*, \ldots, \gamma_4^*$. The periods $f^i(z)$ are a basis
of solutions of the Picard-Fuchs ODE.

Let $\varphi\colon \widetilde{\mathfrak X}\to B$ be a family of
complex structures on $\widetilde X$ with a not necessarily simply
connected base. The fiberwise sublattice $R^3 \varphi_*
{\mathbb{Z}}$ produces a flat connection $\nabla$, the Gauss-Manin
connection, on $R^3\varphi_*{\mathbb C}\otimes{\mathcal O}_B$
which cannot (usually) be globally trivialized. A ${\mathcal
C}^\infty$ trivialization of the family along paths in $B$, with
basepoint $b$, defines a monodromy representation
$$\rho : \pi_1(B,b) \rightarrow \mbox{Aut} \left( H^3(\widetilde{X}_b,
{\mathbb{Z}}), \langle \cdot , \cdot \rangle \right) =: \Gamma \ .$$
Of course $\Gamma$ acts naturally on ${\mathcal D}$, and the
distribution associated to the transversality condition is
$\Gamma$-invariant. The local period maps globalize to a horizontal
holomorphic map
$$\Pi \colon B \rightarrow {\mathcal D} / \Gamma \ .$$
These maps are classifying  Integral Variations of Hodge Structure
(${\mathbb{Z}}$-VHS), i.e., locally liftable horizontal variations
of Hodge structure together with the integral local system, see
\cite[page 165]{Grif}. The image of such maps is contained in the smooth locus
of ${\mathcal D}/\Gamma$. These objects were abstracted and
axiomatized by Deligne \cite{Del1}. Deligne even established a
(non-effective) finiteness theorem for ${\mathbb{Z}}$-VHS over a
base with a fixed discriminant (e.g., over ${\mathbb{P}}^1 \setminus
\{0,1,\infty\}$) \cite{Del3}. In particular Deligne's result implies
that there are finitely many complete horizonal locally liftable
integral curves isomorphic to ${\mathbb{P}}^1 \setminus
\{0,1,\infty\}$ in smooth locus of the period domain. (It is
elementary to see \cite[Theorem (9.8) on page 167]{Grif}
 that there are no such curves isomorphic to
${\mathbb{P}}^1\setminus A$ where $A$ has cardinality at most two.)
Given a holomorphic section $\varpi(z)$ of the line bundle
${\mathcal F}^3$ one still has the Picard-Fuchs equation satisfied
by $\varpi(z)$  and the $\nabla^i\varpi(z), 1\le i\le 4$. There is a
meromorphic map from the trivial bundle of solutions of the
corresponding Picard-Fuchs ODE to the vector bundle given by the
dual local system. This map is a holomorphic isomorphism at every
regular point for the Picard-Fuchs ODE, and it identifies the
monodromy of the Picard-Fuchs ODE with the monomdromy of the
Gauss-Manin connection on the dual local system.

\subsection{Toric examples}

The proposals for just what constitutes a mirror pair of
Calabi-Yau threefolds have evolved significantly since \cite{GP}.
The most commonly used is the Batyrev-Borisov construction
\cite{Bat1,Bat2,Bor,BB1,BB2} in the toric setting, which reduces
the description of mirror pairs of Calabi-Yau hypersurfaces and
complete-intersections in Gorenstein Fano toric varieties to the
classical ``polar duality'' of reflexive polytopes (with NEF
partitions in the complete intersection case).  In this setting
Batyrev-Borisov  \cite{BB2} establish the basic mirror prediction
for Hodge numbers of the proposed mirror pairs, and as a
consequence these pairs are widely believed to be examples of
mirror pairs \cite{MSCMM}.

Let $N$ be a lattice and $M$ the dual lattice. As is well known
\cite{Ful} convex lattice polytopes $\Delta\subset N\otimes
{\mathbb R}$ define toric varieties ${\mathbb{P}}_\Delta =
\mbox{Proj}(S_\Delta)$, where $S_\Delta$ is the polytope ring of
monomials indexed by $n \in k \Delta\cap N$ and graded by total
degree $k$, and where the torus $T_N = N \otimes_{\mathbb Z}
{\mathbb C}^*$ acts. Recall that a lattice polytope $\Delta\subset
N\otimes {\mathbb R}$ is {\em reflexive} if the {\em polar}
polytope
$$\Delta^\circ=\{x\in M\otimes {\mathbb R}\bigl|\bigr. \langle
x,\Delta\rangle\ge -1\}$$ is a lattice polytope with respect to
$M$. When the polytope $\Delta$ is reflexive \cite[\S 3.5]{CoKa},
the toric variety is Fano, and the ``divisor at infinity''
$D_\infty$, i.e., the union of the lower dimensional $T_N$-orbits,
is an anti-canonical divisor. Generic sections of ${\mathcal
O}(D_\infty)$ are then Calabi-Yau threefolds. In general ${\mathbb
P}_\Delta$ is sufficiently singular so that all sections of
${\mathcal O}(D_\infty)$ produce singular threefolds. To remedy
this defect one takes a  simplicial decomposition $\Sigma$ of
$\partial \Delta$ with vertices the set of lattice points. Then
the fan consisting of the cones on the simplices of $\Sigma$
determine a toric variety with isolated singularities\footnote{In
higher dimension these varieties are smooth through codimension
3.} providing a crepant resolution of ${\mathbb P}_{\Delta}$.
Because the singularities are isolated and the resolution is
crepant, the generic section of the pullback of ${\mathcal
O}(D_\infty)$ is a smooth Calabi-Yau $3$-fold. The proposal of
Batyrev \cite{Bat2} is that the Calabi-Yau $X$ obtained from a
maximal triangulation of a reflexive polytope $\Delta$ is mirror
to the Calabi-Yau $\widetilde{X}$ obtained from a maximal
triangulation of the polar $\Delta^\circ$.
\begin{example}\label{reflexives}
The Fano toric variety
${\mathbb{P}}_\Delta = {\mathbb{P}}^4$ is specified by the reflexive
polytope $\Delta$ with vertices
$$\{ (1,0,0,0), (0,1,0,0), (0,0,1,0), (0,0,0,1), (-1,-1,-1,-1) \}
\ .$$  Here the ``divisor at infinity" is a union of the five
coordinate hyperplanes in ${\mathbb{P}}^4$.  The quintic family
${\mathbb{P}}^4[5]$ is interpreted in this context as the family of
hypersurfaces in this (ample) anticanonical divisor class.
The ambient toric variety ${\mathbb{P}}_{\Delta^\circ}$ in which the
quintic mirror family of hypersurfaces are constructed is given by
the polar reflexive polytope $\Delta^\circ$ with vertices
$$\{ (4,-1,-1,-1), (-1,4,-1,-1), (-1,-1,4,-1), (-1,-1,-1,4), (-1,-1,-1,-1) \} \ .$$
This time, however, there are many more integral points in the
polytope besides the vertices and origin and hence one must take a
simplicial decomposition to obtain smooth hypersurfaces. A
combinatorial formula of Batyrev allows one to determine the Hodge
numbers of both the quintic and quintic mirror from the arrangement
of these integral points on the facets of $\Delta$ and
$\Delta^\circ$.

The quintic twin family also has such a toric hypersurface
representation, this time with reflexive polytope the simplex with
vertices
$$\{ (1,0,0,0), (0,1,0,0), (0,0,1,0), (1,2,3,5), (-2,-3,-4,-5) \}
$$
and its polar polytope is
$$\{ (-1,-1,4,-2), (4,-1,-1,0), (-1,4,-1,-1), (-1,-1,-1,2), (-1,-1,-1,1)\} \ .$$
  As before, the families  of both the quintic twin and quintic twin mirror hypersurfaces
  are the
families of  generic sections of the anti-canonical line bundle
over the toric varieties (appropriately resolved). Notice that the
sublattice spanned by the integral points of the quintic twin
lattice polytope has index $5$ in the integral lattice.
\end{example}

What of the associated Picard-Fuchs equations for the family of
holomorphic $3$-forms? In this toric setting the machinery of
Gel'fand-Graev-Zelevinsky \cite{GGZ} and
Gel'fand-Kapranov-Zelevinsky \cite{GKZ} applies to construct, from
the reflexive polytope corresponding to a given Calabi-Yau family, a
differential equation, or more generally a system of such equations
in the case of multiparameter families, satisfied by the periods of
the mirror family (GKZ-system). This is a system of generalized
hypergeometric equations. This system may possess more solutions
than the genuine Picard-Fuchs system, but at least in the
hypersurface case there are methods for describing the reduction to
the actual Picard-Fuchs system \cite{HLY1,HLY2}.

The situation is easier to describe for one-dimensional families
coming from toric geometry. When the family $\widetilde{\mathfrak
X}$ of Calabi-Yau threefolds that are complete intersections in a
Gorenstein Fano toric variety have $h^{2,1}(\widetilde{X}_z)=1$,
then the corresponding Picard-Fuchs ODE is always a generalized
hypergeometric ordinary differential equation \cite[\S 5.5]{CoKa}.
These equations have base ${\mathbb P}^1$ and have  exactly three
regular singular points. The parameter space $B$ of the variation
can always be taken to be ${\mathbb P}^1\setminus\{0,1,\infty\}$.
The local monodromy about the point at infinity in $B$
(conventionally taken to be $z=0$) corresponding to the toric
divisor at infinity $D_\infty$ is maximal unipotent. Let
$\gamma_0$ be an invariant cycle for this monodromy. The period on
$\gamma_0$ of a family
 of holomorphic $3$-forms (defined over all of $B$) is always
represented in the unit disk as a GKZ-generalized hypergeometric
series. A second regular singular point (which we locate at $z=1$)
has local monodromy unipotent of rank one and is called a {\em
conifold singularity} in the physics literature. Of course, by
general principles \cite{BN}, the last has quasi-unipotent
monodromy.

For the quintic mirror and the quintic twin mirror this
generalized hypergeometric series is the one given in
Section~\ref{firstegs}.

\subsection{Homological Mirror Symmetry}\label{HMS}

Much more than a correspondence of Hodge numbers of mirror pairs
is predicted by mirror symmetry.  Na\"ively, for a mirror pair
$(X,\widetilde X)$  deformation of the complex structure on
$\widetilde X$ is mirror to deformation of the symplectic (i.e.,
K\"ahler) structure on $X$. So far there has been no general
construction of a moduli space of symplectic structure
deformations mirror to the usual moduli space of complex structure
deformations. (For such a proposed construction in the case of
$K3$ surfaces, see \cite{Bri}.) In particular, there is no known
analogue for $X$ of the ${\mathbb Z}$-VHS for $\widetilde X$.
Nevertheless, there is a proposal due to Kontsevich which
significantly refines the statement that $b^{\rm even}(X)=b^{\rm
odd}(\widetilde X)$. This is his  Homological Mirror Symmetry
(HMS) Conjecture. It posits that mirror pairs $(X,\widetilde X)$
should have the property that the derived category of bounded
complexes of coherent sheaves on $X$ (clearly a generalization of
the even cohomology) should be equivalent to a category built from
the Fukaya symplectic category of the mirror $\widetilde{X}$,
taking into account certain flat line bundles on the special
Lagrangian cycles (which is a generalization of the third
cohomology of $\widetilde X$). It is not yet a precise conjecture,
and much work in algebraic and symplectic geometry is underway to
make it so. The conjecture was inspired by the earliest inklings
of what eventually became $D$-brane physics, and, through the
introduction of the derived category, has itself provided a
language for physicists to express that theory \cite{MSCMM}. In
the context of the Batyrev-Borisov mirror pairs, illustrative
examples have been studied by many (see for example
\cite{Asp1,Asp2,BD,BDM,DJP,Doug,Hor,Hos1,Hos2,Mor} and references
therein).

As we noted before, there is more structure underlying the family
${\widetilde X}$ than just the Gauss-Manin connection, or
equivalently the Picard-Fuchs differential equation satisfied by
the family of holomorphic $3$-forms. There is the ${\mathbb
Z}$-VHS where the integral local system arises by taking periods
over integral $3$-cycles. This is one of the main objects of our
interest in this investigation.  There is a general result, due to
Griffiths, Deligne, and Schmid, (see the end of this section) that
tells us in this case that the integral monodromy representation
determines the ${\mathbb Z}$-VHS. Now by considering ``$D$-brane
charges'', i.e., by taking cohomology in the categories on both
sides of the HMS correspondence, and retaining the information
about categorical automorphisms on both sides, we obtain a
mathematical prediction from mirror symmetry. By this reasoning,
mirror symmetry predicts that there is an equivalence between the
monodromy representation of the  ${\mathbb Z}$-VHS for
$\widetilde{\mathfrak X}$ and automorphisms of integral even
topological $K$-theory $K^0(X)$ obtained from the categorical
autoequivalences. There is a precise conjecture \cite{Kon} for the
$K$-theory automorphism mirror to the maximal unipotent monodromy
for any Calabi-Yau family. There is also a conjecture
\cite{BD,BDM,DJP} (generalizing a proposal of Kontsevich in the
case of the quintic) for the $K$-theory automorphism mirror to the
monodromy about $z=1$ for one-parameter families, or more
generally for the monodromy about the conifold locus in the
compactification of the moduli space.

 Let us give more details about
these proposed $K$-theory automorphisms. Let $(X,\widetilde X)$ be
a mirror pair with $h^{1,1}(X)=h^{2,1}(\widetilde X)=1$. The
maximal unipotent monodromy for the one-parameter ${\mathbb
Z}$-VHS for $\widetilde{\mathfrak X}$  is proposed to be mirror to
the $K$-theory automorphism $K(X)\to K(X)$ given by
\begin{equation}\label{T0Ktheory} \xi \mapsto \xi \otimes
{\mathcal L}\end{equation}
 where $c_1({\mathcal
L})$ is the positive generator of $H^{1,1}(X,{\mathbb Z})$.  The
conifold monodromy (conventionally taken about $z=1$) is proposed
to be mirror to the Fourier-Mukai ``push-pull'' transform
\begin{equation}\label{T1Ktheory}
\xi \mapsto \sum_\lambda (p_1)_* \left( \{ (\lambda^{-1} \boxtimes
\lambda) \rightarrow {\mathcal O}_\Delta \} \otimes p_2^* (\xi)
\right)\end{equation} where $\Delta \subset X \times X$ is the
diagonal, $p_i$ are the projections onto the factors, and
$\lambda$ runs over flat line bundles on $X$. This proposal
implies that the monodromy automorphism of $H^3(\widetilde{X},
{\mathbb Z})$ around the conifold locus is then a sum of terms
each of which is a Picard-Lefschetz transformation on a vanishing
cycle.  The vanishing cycles are all equivalent modulo torsion, so
that this monodromy is divisible by the number of terms in the
sum, which is $|\pi_1(X)|$.  The remaining monodromy
(conventionally taken at $\infty$), whether of finite order (as in
the case of the quintic family) or not, is determined as the
product of the monodromies at $0$ and $1$.

\section{Classification of ${\mathbb Z}$-VHS}

For the purpose of investigating this proposed
Hodge-theoretic/$K$-theoretic mirror relationship it is very useful
to have a complete description of the possible integral variations
of Hodge structure that can underlie a family of Calabi-Yau
threefolds over ${\mathbb{P}}^1 \setminus\{0,1,\infty\}$ with
$h^{2,1} = 1$, subject to the assumption of the existence of a point
of maximal unipotent monodromy (at $z=0$), a point of unipotent
monodromy of rank one (at $z=1$), and quasi-unipotency of their
product (necessary at $z=\infty$ by the Monodromy Theorem \cite{BN}.
In this case  the monodromy representation of the ${\mathbb Z}$-VHS
determines that variation up to isomorphism. To see this let
${\mathcal H}$ and ${\mathcal H}'$ be ${\mathbb Z}$-VHS over the
same connected algebraic base $B$ with $h^{3,0}=h^{2,1}=1$ and with
isomorphic irreducible monodromy ((2) of Theorem~\ref{real}
establishes irreducibility for our specific case). The bundle
${\mathcal Hom}({\mathcal H},{\mathcal H}')$ with its natural
induced flat connection, integral structure and Hodge filtration  is
a ${\mathbb Z}$-VHS of weight zero over $B$. By \cite[Theorem
(7.22)]{Sch}, the $\pi_1(B)$-invariant subbundle is a sub ${\mathbb
Z}$-VHS. Of course, the invariant subbundle is the flat bundle whose
fiber over any point $x\in B$  is $\pi_1(B,x)$-invariant subspace of
$\Hom_{\mathbb C}({\mathcal H}_x,{\mathcal H}'_x)$. Since the
monodromies of ${\mathcal H}$ and ${\mathcal H}'$ are irreducible ,
and hence by Shur's lemma, the invariant subbundle is
one-dimensional. But a weight zero Hodge structure on a
one-dimensional space is of type $(0,0)$. This means that any
non-trivial $\pi_1(B,x)$-invariant homomorphism ${\mathcal H}_x\to
{\mathcal H}'_x$ preserves the Hodge filtration and hence is an
isomorphism of Hodge structures, showing that ${\mathcal H}$ and
${\mathcal H}'$ are isomorphic as ${\mathbb Z}$-VHS over $B$.

In this section we classify ${\mathbb Z}$-VHS over ${\mathbb
P}^1\setminus\{0,1,\infty\}$ with $h^{3,0}=h^{2,1}=1$ subject to the
conditions that the local monodromy about $z=0$ is maximal unipotent
and the local monodromy about $z=1$ is unipotent of rank one. As
noted above, the key to this classification is the classification of
possible monodromy representations underlying such variations.

\subsection{Classification results for $\rho\colon
\pi_1({\mathbb{P}}^1 \setminus \{0,1,\infty\})\to
Sp(4,{\mathbb{Z}})$}\label{reps}

We begin by classifying real representations. Let $V$ be a four
dimensional real vector space. Let $T\colon V\to V$ be a unipotent
automorphism of  $V$ and let $N$ be the nilpotent endomorphism of
$V$ defined by $N=T-\Id$. Then $N^{{\rm dim}V}=0$. $T$ is said to be
{\sl maximal unipotent} if $N^{{\rm dim}V-1}\not=0$. Also, the {\sl
rank} of $T$ is defined to be the dimension of the image of $N$, so
that $T$ is maximal unipotent if and only if its rank is one less
than the dimension of $V$.

Throughout we fix $x\in {\mathbb{P}}^1 \setminus \{0,1,\infty\}$ and
a representation
$$\rho\colon \pi_1({\mathbb{P}}^1 \setminus
\{0,1,\infty\},x)\to GL(V).$$ We have geometric generators
$\gamma_0$,$\gamma_1$, and $\gamma_\infty$ for $\pi_1(X,x)$  which
are represented by loops in ${\mathbb{P}}^1 \setminus
\{0,1,\infty\}$ encircling $0,1,\infty$ respectively. These are
oriented positively  so that they satisfy the relation
$\gamma_0\gamma_1\gamma_\infty=1$. We denote the images of these
generators under $\rho$ by $T_0,T_1,T_\infty$, respectively, and in
view of the eventual application to flat bundles over
${\mathbb{P}}^1 \setminus \{0,1,\infty\}$ call them the {\sl
monodromy} around $0,1,\infty$.

\begin{theorem}\label{real}
Suppose that the monodromy $T_0$  is maximal unipotent and that the
monodromy $T_1$  is unipotent of rank $1$. Set $N_i=T_i- Id$ for
$i=0,1$.  Define an invariant $m\in \Ar$ by choosing a non-zero
vector $v\in \Ker N_0$ and setting $N_0^3(N_1(v))=-mv.$
 Then:
 \begin{enumerate}
 \item $m$ is independent of the choice of $v$.
 \item The representation $\rho$ is irreducible if and only if
 $m\not= 0$.
 \item If the  monodromy at infinity has an
 invariant vector, then $m=0$.
 \item In the case when $m\not= 0$, the
 representation  $\rho$ is determined, up to conjugation by an element
 of $GL(V)$, by the
 characteristic polynomial of the monodromy transformation around
 infinity.  Direct computation shows that the characteristic polynomial of
$T_\infty^{-1}=T_0T_1$ is
\begin{equation}\label{charpoly}
x^4+(a-4)\,x^3+(6-2a+bm)\,x^2+(a-4+m-bm)\,x+1 \ . \end{equation}
Notice that $a$ and $b$ are determined by this polynomial.  The real
numbers $a,b,m$ are complete invariants of the conjugacy class of
$\rho$ as a representation into $GL(V)$.
 \item For $m\not=0$
 take as basis for $V$ the vectors
\begin{equation}
N_1(-v), N_0(N_1(-v)), \frac{N_0^2(N_1(-v))}{m},v \ .
\label{Vbasis}
\end{equation}
Then in this basis matrices for $T_0$ and $T_1$ are:
 \begin{eqnarray}\label{1}
T_0&=&\begin{pmatrix}1&0&0&0\\1&1&0&0\\0&m&1&0\\0&0&1&1\end{pmatrix}
\\ T_1 & = &
\begin{pmatrix}1&-a&-b&-1\\0&1&0&0\\0&0&1&0\\0&0&0&1\end{pmatrix}\label{2}
\end{eqnarray}
 \end{enumerate}
 \end{theorem}

{\bf For the rest of this section we assume that $T_0$ is maximal
unipotent, that $T_1$ is a unipotent of rank one, and that $\rho$ is
irreducible, i.e., $m\not=0$, and we let $m,a,b$ be as in the
statement of the previous lemma.}

\begin{corollary}\label{b=1}
There is a non-degenerate symplectic form on $V$ invariant under
$\rho$ if and only if $b=1$. In this case, this form is unique up
to non-zero scalar multiplication and the matrix for this form in
the basis (\ref{Vbasis}) is (up to a positive scalar):
\begin{equation} \label{3}
\langle
\cdot,\cdot\rangle= \pm \begin{pmatrix} 0&-a&-1&-1\\a&0&1&0\\1&-1&0&0\\1&0&0&0
\end{pmatrix}
\end{equation}
\end{corollary}

{\bf We now fix a non-degenerate skew pairing
$\langle\cdot,\cdot\rangle$ on $V$.}

\begin{corollary}\label{the14}
There are exactly 14 conjugacy classes of representations
$\rho\colon \pi_1({\mathbb P}^1 \setminus \{0,1,\infty\})\to Sp(V)$
with $T_0$ and $T_1$ as given in  Theorem~\ref{real} and with
$T_\infty$ being quasi-unipotent. These are classified by the pairs
of nonzero integers $(m,a)$ (with $b=1$) and a complete set of such
pairs is given in the  first two columns of Table \ref{thetable}.
\end{corollary}

We are actually interested in representations to $Sp(4,{\mathbb
Z})$, or equivalently in lattices $L\subset V$ invariant under
$\rho$ and on which the symplectic form is a perfect integral
pairing. Here is a first result along these lines.

\begin{corollary}\label{intexist}
 Suppose that $L\subset V$ is
an integral lattice invariant under $\rho$ such that the
restriction of  $\langle \cdot,\cdot \rangle$ to $L$ is a perfect
integral pairing. Then  $b=1$ and $a,m\in \Zee$. Conversely, if
$b=1$ and $a,m\in \Zee$, then  there is an integral lattice $L_0$
in $V$ invariant under $\rho$
 on which the pairing is a perfect integral pairing. Furthermore, there
 is a basis for $L_0$ in which
  $T_0$, $T_1$ and
 the pairing   are given by the matrices in
 Equations~(\ref{1}),and~(\ref{2}), respectively, and the symplectic form is
given,
 up to sign, by Equation~(\ref{3}).
\end{corollary}

Notice that we are not claiming that given a lattice $L\subset V$
invariant under $\rho$ and on which the symplectic form is a perfect
integral pairing, the lattice $L_0$ as in Corollary~\ref{intexist}
is related in any way to $L$. In fact, they are closely related
(though not necessarily equal), but that relationship is somewhat
delicate and depends on the arithmetic properties of $m$ and $a$.

Now suppose that $b=1$ and that $a$ and $m\not= 0$ are integers. Let
$e_1,e_2,e_3,e_4$ be a basis of $V$ in which the symplectic pairing
$\langle \cdot, \cdot \rangle$ is given by Equation~(\ref{3}), and
$T_0$ and $T_1$  are given by Equations~(\ref{1}) and~(\ref{2}),
respectively. Our goal is to study all possible integral lattices
invariant under $\rho$ and on which the symplectic pairing is a
perfect integral pairing. By Corollary~\ref{intexist} the lattice
$L_0$ spanned by $\{e_1,\ldots,e_4\}$ is such a lattice.

We introduce the (increasing) weight filtration $W_*(V)$ associated
with $T_0$ by setting, $W_{2i}=W_{2i+1}=\Ker(N_0)^{i+1}$ for
$i=0,\ldots,3$. This is the monodromy weight filtration in the
limiting mixed Hodge structure as defined by Deligne \cite{Del1}.
This filtration has the property that $N_0\colon W_j\to W_{j-2}$ and
$N_0$ induces an isomorphism on the associated gradeds: $N_0\colon
W_j/W_{j-1}\to W_{j-2}/W_{j-3}$. Fix an integral lattice $L$ on
which the pairing is a perfect integral pairing and which is
invariant under $\rho$. Then $W_{2i}(L)=L\cap W_{2i}$ has rank equal
to $i+1$, and hence the quotients $W_{2i}(L)/W_{2i-2}(L)\cong \Zee$
for $i=0\ldots,3$.
 Notice in particular that for the lattice $L_0$
$$N_1\colon W_0(L_0)\to L_0/W_4(L_0)$$  and
$$N_0\colon L_0/W_4(L_0)\to W_4(L_0)/W_2(L_0)$$
are isomorphisms.

\begin{lemma}\label{lowertrian}
Suppose  that $\ell_1,\ell_2,\ell_3,\ell_4$ is a basis for $L$
such that for each $i$ the subset $\ell_{4-i},\ldots,\ell_4$ a
basis for $W_{2i}(L)$. Then the matrix for the symplectic form
is: $$\begin{pmatrix} 0 &  \alpha & \beta & \pm 1 \\
-\alpha & 0 &\pm 1 & 0 \\ -\beta & \pm 1 & 0 & 0 \\
\pm 1 & 0 & 0 & 0\end{pmatrix}$$ for appropriate integers $\alpha$
and $\beta$.\end{lemma}

\begin{remark} We now fix our  conventions regarding integral bases for $L$
pertaining to the above isomorphisms on associated gradeds. From now
on we shall consider bases for $L$ satisfying: \begin{enumerate}
\item For every $i\le 4$ the set $\{\ell_{4-i},\ldots,\ell_4\}$
forms a basis for $W_{2i}(L)$; \item $\langle \ell_1, \ell_4 \rangle
=-1$; \item $\langle \ell_2, \ell_3 \rangle =1$. \item Duality tells
us that there is a non-zero integer $r$ such that $N_0(\ell_1)\equiv
r\ell_2 \pmod{\ell_3,\ell_4}$ and $N_0(\ell_3)=r\ell_4$. By changing
the signs of $\ell_2$ and $\ell_3$ if necessary, we arrange that
$r>0$. \item  We define $s$ by requiring that $N_0(\ell_2)\equiv
s\ell_3 \pmod{\ell_4}$. \item We define $t$ such that
$N_1(\ell_4)=-t\ell_1\pmod{\ell_2,\ell_3,\ell_4}$.
\end{enumerate}
Of course, $r$ is  the divisibility of $N_0$ as a map $W_6/W_4\to
W_4/W_2$ and also $W_2/W_0\to W_0$ while $s$ is the divisibility of
$N_0\colon W_4/W_2\to W_2/W_0$ and $t$ is the divisibility of
$N_1\colon W_0\to W_6/W_4$.
\end{remark}

In particular, we have:

\begin{theorem}
The integers $r,s$ and $t$ are  invariants of the isomorphism
class of $(L,\rho,\langle \cdot, \cdot \rangle)$.
\end{theorem}

Notice that these integral invariants are related to the integer
invariant $m$ of the real conjugacy class of $\rho$ by:
\begin{equation}\label{relation}
m=r^2st \ \ \ {\rm and}\ \ \  t|a.
\end{equation}

For each $i$, $1\le i\le 4$, the real subspace of $V$ spanned by
$\{\ell_i,\ldots \ell_4\}$ is equal to that spanned by
$\{e_i,\ldots,e_4\}$, which means the matrix whose columns express
the $\ell_i$ in terms of the $e_j$ is weakly lower triangular.
Since the image of the representation $\rho$ is Zariski dense in
$Sp(V)$
 general considerations imply that the coefficients in this matrix are
 contained in $\Q[\sqrt{(1/d)}]$ for some integer $d$.
 In fact, a much more explicit result holds.

As an indication of how the basis $\{\ell_i\}$ of the new lattice
$L$ is expressed in terms of the original basis $\{e_i\}$ of the
original lattice $L_0$, consider the following lemma.
\begin{lemma}\label{root} There exists a basis
$\ell_1,\ell_2,\ell_3,\ell_4$ for $L$ as described above such that
\begin{eqnarray*} \ell_1 & = & \frac{e_1}{\sqrt{t}}+\frac{\alpha
e_2}{rt\sqrt{t}}+ \frac{\beta e_3}{\sqrt{t}}+\frac{\gamma e_4}{rt\sqrt{t}}
\\
\ell_2 & = & \frac{e_2}{r\sqrt{t}}+\frac{\delta e_3}{\sqrt{t}}+
\frac{\mu e_4}{r\sqrt{t}}\\
\ell_3 & = & r\sqrt{t}e_3  +\frac{\alpha e_4}{\sqrt{t}} \\
\ell_4 & = & \sqrt{t} e_4
\end{eqnarray*}
for appropriate $\alpha,\beta,\gamma,\delta,\mu\in \Zee$.
\end{lemma}

Once again let us impose the condition that $T_\infty$ is
quasi-unipotent, i.e., that the pair $(m,a)$ is one of the 14 pairs
listed in Table~\ref{thetable}. Case-by-case analysis results in a
complete classification of the possible $\alpha, \beta, \gamma,
\delta, \mu$ corresponding to each allowable $r,s,t$. Adding the
numbers in Column 4 of Table~\ref{thetable} one see that there are
112 possibilities.

In Equations~(\ref{T0Ktheory}) and~(\ref{T1Ktheory}) we gave
explicit formulas for the proposed automorphisms of $K$-theory.
Translating these by mirror symmetry gives the following
conditions on $T_0$ and $T_1$, which we express in terms of the
invariants $r$ and $t$:
\begin{conj}\label{conjecture}
For any ${\mathbb Z}$-VHS arising from a family of Calabi-Yau
threefolds $\widetilde{\mathfrak X}$ with $h^{2,1}=1$ over
${\mathbb P}^1\setminus\{0,1,\infty\}$ we have:
\begin{eqnarray}
r & = & 1 {\rm \ \ (mirror\ to\ the\ fact\ that\ } c_1({\mathcal
L})
{\rm \ is\ a\ generator\ of\ } H^2(X,{\mathbb Z}) {\rm )\ }  \label{item1} \\
\label{item2}  N_1 & = & t\cdot \overline N_1 \, {\rm \ for\
some\ indivisible\ integral\ transformation\ } \overline N_1.
\end{eqnarray}
Furthermore, $t$  is equal to the order of the fundamental group
of the mirror Calabi-Yau threefold $X$.
\end{conj}

\begin{remark}
A stronger version of the first equation  was conjectured by
Morrison to hold for the local maximal unipotent monodromy in
arbitrary dimensional families of Calabi-Yau threefolds,
\cite[p.~150]{CoKa}.  Condition~(\ref{item2}) and the last
statement are conjectured to hold for the local monodromy around
the principal discriminant for arbitrary dimensional families of
Calabi-Yau threefolds.

The invariants $m$ and $a$ encode natural geometric information
about the Calabi-Yau threefolds $X$ in the last column of
Table~\ref{thetable}.  As observed by Borcea \cite{Borc} for the
13 examples with $t=1$, the integers we identified as invariants
$(m,a)$ above have a mirror interpretation as $m = {\mathcal L}^3$
and $a = \dim(H^0(X,{\mathcal O}({\mathcal L})))$ respectively. In
fact, this interpretation applies to all the geometric examples
listed in Table~\ref{thetable}.
\end{remark}

The classification of representations into $Sp(4,{\mathbb Z})$
satisfying the two conditions in Conjecture \ref{conjecture}  is
simpler to state than the general case:
\begin{prop}
An integral lattice $L\subset V$ on which the symplectic pairing is
a perfect integral pairing and which is invariant under $\rho$ and
which satisfies Conditions~(\ref{item1}) and ~(\ref{item2}) is
determined by $t$.
\end{prop}

Columns 5 and 6 of Table~\ref{thetable} list the 23 of the 112
representations that satisfy these conditions. Notice that without
assuming Conditions~(\ref{item1}) and~(\ref{item2}) the invariants
$m,a,r,s,t$ do not in general determine the isomorphism class of the
lattice.

\subsection{From representations to ${\mathbb Z}$-VHS}

As we argued at the end of Section~\ref{HMS}, for a ${\mathbb
Z}$-VHS over a quasi-projective base with completely irreducible
monodromy representation, the isomorphism class of the monodromy
representation determines the isomorphism class of the ${\mathbb
Z}$-VHS. According to a theorem of Borel and Narasimhan \cite{BN}
the local monodromy of any ${\mathbb Z}$-VHS is quasi-unipotent.

\begin{theorem}\label{thmreduce}
Suppose we have a primitive weight-three  ${\mathbb Z}$-VHS
${\mathcal H}$ over ${\mathbb P}^1\setminus\{0,1,\infty\}$
satisfying:
\begin{enumerate}
\item $h^{3,0}=h^{2,1}=1$,
\item the underlying monodromy representation $\rho$ into $Sp(4,{\mathbb Z})$
is irreducible,
\item $T_0$ is maximal unipotent, and
\item $T_1$ is unipotent of rank one.
\end{enumerate}
Then $\rho$ is conjugate in $Sp(4,{\mathbb Z})$ to one of the 112
representations enumerated in Table~\ref{thetable} and the
isomorphism class of ${\mathcal H}$ as a ${\mathbb Z}$-VHS is
determined by the $Sp(4,{\mathbb Z})$-conjugacy class of $\rho$.
\end{theorem}

The remaining question is one of existence: Which of the 112
representations underlie a ${\mathbb Z}$-VHS?

Each of the 14 real representations comes from a generalized
hypergeometric ODE. That is to say for each of the fourteen pairs
$(m,a)$ in Table~\ref{thetable} we have the following data:
\begin{enumerate}
\item a flat bundle with a real structure,
\item a real symplectic form parallel under the flat connection,
\item a holomorphic section $\varpi$ of the complexification of the flat bundle,
\end{enumerate}
such that the Picard-Fuchs equation satisfied by $\varpi$ is
generalized hypergeometric equation:
\begin{equation}\label{hypgeoeqn}
\left[ \Theta^4 - z \left( \Theta + a_1 \right) \left( \Theta +
a_2\right) \left( \Theta + a_3 \right) \left( \Theta + a_4 \right)
\right] \varpi(z) = 0 \ ,\end{equation}
 where the $a_i$ are the
entries in Column 3 of Table~\ref{thetable} corresponding to $(m,a)$
and $\Theta=z\nabla_z$. The putative Hodge filtration on this flat
bundle is given by: ${\mathcal F}^{3-i}$ is the span of the
holomorphic section and its first $i$ covariant derivatives. This
filtration is automatically horizontal and satisfies the first Hodge
condition~(\ref{hodge1}). It is not clear that the filtration
satisfies the second Hodge condition~(\ref{hodge2}), but if it does,
then the polarization condition~(\ref{hodge3}) follows.

 By direct geometric constructions we know
that all 14 of these representations arise as local systems (or
subquotients of local systems) underlying families of Calabi-Yau
threefolds.  Of these, 13 are local systems of third cohomology
for one-parameter families of complete intersection Calabi-Yau
threefolds with $h^{2,1}=1$ in toric varieties (see Section
\ref{geomreal}). Thus, these thirteen are in fact real
representations underlying ${\mathbb R}$-VHS. Consequently, for
these 13 real types there is a Hodge filtration making this a
${\mathbb R}$-VHS. In the remaining case for which $(m,a) = (1,4)$
we do not know whether Conditions~(\ref{hodge2})
and~(\ref{hodge3}) hold.

Since Conditions~(\ref{hodge1}), (\ref{hodge2}), and~(\ref{hodge3})
make reference only to the symplectic form and the real structure, a
parallel lattice in a ${\mathbb R}$-VHS on which the symplectic form
is a perfect integral pairing determines a ${\mathbb Z}$-VHS. Hence,
all integral representations whose real types satisfy $(m,a)\not=
(1,4)$ underlie (unique isomorphism classes of) ${\mathbb Z}$-VHS.
Since the  real representation corresponding to $(m,a)=(1,4)$ has a
unique integral lattice, there is only one of the 112 integral
representations enumerated in Table~\ref{thetable} which might not
come from a ${\mathbb Z}$-VHS. Thus, up to this one ambiguity we
have completely classified ${\mathbb Z}$-VHS satisfying the
conditions given in Theorem~\ref{thmreduce}. This means that, up to
the same ambiguity, we have classified the complete curves in the
smooth locus of the period space ${\mathcal D}/\Gamma$ isomorphic to
${\mathbb P}^1\setminus \{0,1,\infty\}$ with the given local
monodromy conditions.

\section{Geometric realization of abstract variations by
one-parameter families of Calabi-Yau threefolds} \label{geomreal}

\begin{table}
\begin{center}
\begin{tabular}{||c|c|c|c|c|c|l||}
\hline \rule{0pt}{5mm} $m$ & $a$ & $a_1, a_2, a_3, a_4$
& \# $L_{\mathbb Z}$ & \# $L_{\rm MC}$ & $t$ & Geometric examples \\[3pt]
\hline \hline \rule{0pt}{5mm} 1 & 4 & $\frac{1}{12}, \frac{5}{12},
\frac{7}{12}, \frac{11}{12}$
& 1 & 1 & 1 & I \\[3pt]
\hline \rule{0pt}{5mm} 1 & 3 & $\frac{1}{10}, \frac{3}{10},
\frac{7}{10}, \frac{9}{10}$
& 1 & 1 & 1 & ${\mathbb{WP}}^4_{1,1,1,2,5}[10]$ \\[3pt]
\hline \rule{0pt}{5mm} 2 & 4 & $\frac{1}{8}, \frac{3}{8},
\frac{5}{8}, \frac{7}{8}$ & 2 & 2 & 1 &
${\mathbb{WP}}^4_{1,1,1,1,4}[8]$ \\[3pt]
\cline{6-7} \rule{0pt}{5mm}
 & &  & & & 2 &  II  \\[3pt]
\hline \rule{0pt}{5mm} 5 & 5 & $\frac{1}{5}, \frac{2}{5},
\frac{3}{5}, \frac{4}{5}$ & 2 &
2 & 1 & ${\mathbb{P}}^4[5]$ \\[3pt]
\cline{6-7} \rule{0pt}{5mm}
 & & & & & 5 & A \\[3pt]
\hline \rule{0pt}{5mm} 1 & 2 & $\frac{1}{6}, \frac{1}{6},
\frac{5}{6}, \frac{5}{6}$ & 1 & 1 & 1 &
${\mathbb{WP}}^5_{1,1,2,2,3,3}[6,6]$ \\[3pt]
\hline \rule{0pt}{5mm} 2 & 3 & $\frac{1}{6}, \frac{1}{4},
\frac{3}{4}, \frac{5}{6}$ & 1 & 1 & 1 &
${\mathbb{WP}}^5_{1,1,1,2,2,3}[4,6]$ * \\[3pt]
\hline \rule{0pt}{5mm} 3 & 4 & $\frac{1}{6}, \frac{1}{3},
\frac{2}{3}, \frac{5}{6}$ & 1 & 1 & 1 &
${\mathbb{WP}}^4_{1,1,1,1,2}[6]$ \\[3pt]
\hline \rule{0pt}{5mm} 4 & 5 & $\frac{1}{6}, \frac{1}{2},
\frac{1}{2}, \frac{5}{6}$ & 11 & 1 &
1 & ${\mathbb{WP}}^5_{1,1,1,1,1,3}[2,6]$ \\[3pt]
\hline \rule{0pt}{5mm} 4 & 4 & $\frac{1}{4}, \frac{1}{4},
\frac{3}{4}, \frac{3}{4}$ & 8 & 3 & 1 &
${\mathbb{WP}}^5_{1,1,1,1,2,2}[4,4]$ \\[3pt]
\cline{6-7} \rule{0pt}{5mm}
 & & & & & 2 & B \\[3pt]
\cline{6-7} \rule{0pt}{5mm}
 & & & & & 4 & C \\[3pt]
 \hline \rule{0pt}{5mm}
6 & 5 & $\frac{1}{4}, \frac{1}{3}, \frac{2}{3}, \frac{3}{4}$ & 1 &
1 & 1 & ${\mathbb{WP}}^5_{1,1,1,1,1,2}[3,4]$ * \\[3pt]
\hline \rule{0pt}{5mm} 8 & 6 & $\frac{1}{4}, \frac{1}{2},
\frac{1}{2}, \frac{3}{4}$ & 14 & 2 & 1 &
${\mathbb{P}}^5[2,4]$ \\[3pt]
\cline{6-7} \rule{0pt}{5mm}
 & & & & & 2 & D \\[3pt]
\hline \rule{0pt}{5mm} 9 & 6 & $\frac{1}{3}, \frac{1}{3},
\frac{2}{3}, \frac{2}{3}$ & 8 & 2 & 1 &
${\mathbb{P}}^5[3,3]$ \\[3pt]
\cline{6-7} \rule{0pt}{5mm}
 & & & & & 3 & E \\[3pt]
 \hline \rule{0pt}{5mm} 12 & 7 & $\frac{1}{3}, \frac{1}{2},
 \frac{1}{2}, \frac{2}{3}$ & 11 & 1 & 1 & ${\mathbb{P}}^6[2,2,3]$
 \\[3pt]
 \hline \rule{0pt}{5mm}
 16 & 8 & $\frac{1}{2}, \frac{1}{2}, \frac{1}{2}, \frac{1}{2}$  & 50
 & 4 & 1 &
${\mathbb{P}}^7[2,2,2,2]$ \\[3pt]
\cline{6-7} \rule{0pt}{5mm} & & & & & 2 & F \\[3pt]
\cline{6-7} \rule{0pt}{5mm}  & & & & & 4 & G \\[3pt]
 \cline{6-7} \rule{0pt}{5mm} & & & & & 8 & H \\[3pt]
  \hline \hline
\end{tabular}
\end{center} \caption{Table of ``Mirror-Consistent'' $\rho\colon
\pi_1({\mathbb{P}}^1 \setminus \{0, 1, \infty\}) \to Sp(4)$}
\label{thetable}
\end{table}

The first two columns of Table \ref{thetable}  give the invariants
$m$ and $a$ for the $14$ ${\mathbb{R}}$-VHS in the classification.
The third column gives the coefficients of the hypergeometric ODE,
see Equation~(\ref{hypgeoeqn}), or equivalently  $(1/2\pi \imath)$
times the logs of the eigenvalues of $T_\infty$.  The fourth
column enumerates the integral lattices within each real
representation. The fifth column indicates the number of integral
lattice that are ``mirror-consistent'' in the sense that they
satisfy both conditions (\ref{item1}) and (\ref{item2}). The sixth
column shows the $t$ values for each of the mirror consistent
representations enumerated in the fifth. For each of these, the
final column indicates known geometric examples with integral
variations of Hodge structure of this type. Since the examples
listed all have $h^{1,1} = 1$, to read each entry preface with
``The ${\mathbb{Z}}$-VHS of $H^3$ of the Batyrev-Borisov mirror of
$\ldots$''.

The two examples with $(m,a) = (2,3)$ and $(6,5)$ have an asterisk
to indicate that, even though these families of Calabi-Yau
complete intersections in weighted projective space are of the
correct type, the Newton polytopes of the weighted projective
spaces are not reflexive and so these examples don't quite fit
into the Batyrev-Borisov description of mirror pairs (though the
Greene-Plesser construction does apply \cite{KT}). Denoting the
vertices of the $(2,3)$ polytope by $\{e_1, \ldots, e_5,
-e_1-e_2-2e_3-2e_4-3e_5\}$, and those of the $(6,5)$ polytope by
$\{e_1, \ldots, e_5, -e_1-e_2-e_3-e_4-2e_5\}$, suitable reflexive
polytopes are obtained by adding in each case the new vertex
$-e_5$. These have the property that they possess NEF partitions
of bidegrees $[4,6]$ and $[3,4]$ respectively, with Hodge numbers
matching those of the original Calabi-Yau hypersurfaces in the two
weighted projective spaces.

The example denoted  A is the quintic twin, which we already
recognized as the free quotient of the quintic by ${\mathbb
Z}/5{\mathbb Z}$. The rest of the examples labeled B through H arise
in a quite similar fashion, by finding a suitable free action by a
group of order $t \neq 1$ on a simply-connected Calabi-Yau of that
class (specifically the one described in the corresponding $t=1$
entry of the fifth column) and taking the quotient.  When that
action is compatible with the toric structure of the original
example, the Batyrev-Borisov mirror construction applies. In the
classification of 4D reflexive polytopes by Kreuzer-Skarke [KrSkCY]
there are precisely five polytopes which correspond to toric
hypersurfaces with $h^{1,1}=1$, and four of these are among the
weighted projective spaces listed in the seventh column of Table
\ref{thetable}
--- the quintic, the sextic, the octic, and the dectic.  The
remaining example is the alternative derivation of the quintic
twin from Example~\ref{reflexives} for which the Batyrev mirror
pair construction applies. Note that the fact the ``octic twin''
does not arise in this way means that if there does exist a free
quotient of the octic hypersurface in
${\mathbb{WP}}^4_{1,1,1,1,4}[8]$ by the group
${\mathbb{Z}}/2{\mathbb{Z}}$, then the action cannot be compatible
with the toric structure (i.e., arise from a natural involution on
the weighted projective ambient space).

The only known construction of example H, a Calabi-Yau whose
mirror defines the $t=8$ integral structure in the same real class
as the mirror of the complete intersection of four quadrics in
${\mathbb{P}}^7$, requires a manifestly non-toric group action. It
comes from Jae Park's recent construction \cite{Park} of the
mirror of Beauville's Calabi-Yau threefold with nonabelian
fundamental group. Example H is the Beauville manifold, the
quotient of ${\mathbb{P}}^7[2,2,2,2]$ by the group $Q_8$ of unit
quaternions. Examples F and G are then the intermediate quotients
obtained through the same action, but restricted to the subgroups
${\mathbb{Z}}_2$ and ${\mathbb{Z}}_4$ of $Q_8$ respectively.

For Example II the authors don't know of a family of Calabi-Yau
threefolds with $h^{2,1} = 1$ over ${\mathbb{P}}^1 \setminus \{0, 1,
\infty\}$ realizing this ${\mathbb{Z}}$-VHS, and for Example I the
authors don't know whether there is a geometric realization of this
representation or indeed whether there is a ${\mathbb Z}$-VHS with
this representation.

\begin{remark}
All of these toric examples satisfy Conjecture~\ref{conjecture}.
Notice that geometric families realizing a given $(m,a)$ but
different $t$ values give examples of geometric ${\mathbb Z}$-VHS
that are isomorphic as ${\mathbb R}$-VHS.  The quintic mirror and
the quintic twin mirror families are the simplest example of this
phenomenon.  More generally, to the best of the authors'
knowledge, all complete intersection Calabi-Yau threefolds arising
from a NEF partition of the anti-canonical divisor at infinity in
Gorenstein Fano toric varieties satisfy the multiparameter
generalization of this conjecture.
\end{remark}

Some positive results along these lines are:

\begin{theorem}
For any Calabi-Yau threefold that is an anti-canonical
hypersurface in a Gorenstein Fano toric variety the local
monodromy $T_1={\rm Id}+N_1$ around the prinicpal discriminant
satisfies Condition~(\ref{item2}) with $t$ equal to the order of
the fundamental group of the Batyrev mirror Calabi-Yau.

For any Calabi-Yau threefold $\widetilde X$ that is a complete
intersection arising from a NEF partition of the anti-canonical
divisor at infinity in a Gorenstein Fano toric variety and with
$h^{2,1}(\widetilde X)=1$ the parameter space of complex moduli
can be identified with ${\mathbb P}^1\setminus \{0,1,\infty\}$ in
such a way that the monodromies $T_0$ and $T_1$ satisfy the
assumptions of Theorem~\ref{thmreduce}. Furthermore, $T_0$
satisfies Condition~(\ref{item1}).
\end{theorem}

There are five one-parameter families of Calabi-Yau threefold
hypersurfaces in a Gorenstein Fano toric variety with $h^{2,1} =
1$.  By this theorem each of their ${\mathbb Z}$-VHS is one of the
23 listed in column 5 of Table~\ref{thetable}, where $(m,a)$ are
determined by the local monodromy matrices $T_0$ and $T_1$ and $t$
is the order of the fundamental group of the mirror.  We
conjecture that the same is true for Calabi-Yau threefold complete
intersections in Gorenstein Fano toric varieties with $h^{2,1} =
1$.

\section{Extensions and final comments}

\subsection{Extensions}

The classification obtained above for ${\mathbb{Z}}$-VHS underlying
$h^{2,1}=1$ Calabi-Yau threefold moduli over ${\mathbb{P}}^1
\setminus \{0,1,\infty\}$ should be extended both to $h^{2,1} > 1$
multiparameter hypergeometric variations and to cases of $h^{2,1}=1$
Calabi-Yau threefold moduli over curves of higher genus and with
more punctures.

A small but growing number of the latter sorts of Calabi-Yau
threefold families are known to exist, mostly via complete
intersection constructions in partial flag varieties \cite{BCFKvS},
though there is a great deal of indirect evidence for their
existence. In particular, recent work of Almkvist-Zudilin \cite{AZ}
on ordinary differential equations whose associated formal ``Yukawa
coupling'' series possess integrality and positivity properties that
suggest they could be Yukawa couplings for actual Calabi-Yau
threefolds
--- so-called ``fake Picard-Fuchs equations'' --- highlights the
importance of classifying compatible integral structures, both as
a means of eliminating some cases from consideration and in order
to obtain more information about compatible Calabi-Yau manifolds
if they do exist (e.g., possible nontrivial fundamental groups).

For a large number of the Almkvist-Zudilin examples (with $t=1$)
Christian van Enckevort and Duco van Straten \cite{vEvS} have used
a numerical approach, motivated by HMS and the integrality of
formal  Gopakumar-Vafa invariants, to determine the integral
monodromy representations. One very interesting question is
whether their methods can be modified to also detect alternate
integral structures with $t > 1$.  These would correspond
conjecturally to more families of mirrors of non-simply connected
Calabi-Yau threefolds.

For each integral structure with $t > 1$ in a class for which a
geometric construction of the $t = 1$ representative is known, one
can ask whether there is a free action of a finite group of order
$t$ on the known example such that the mirror family realizes the
$t > 1$ variation.  In particular, for the examples in
Table~\ref{thetable} the only entry for which this remains an open
question is case II.  Thus we ask: Is there a free, holomorphic
$3$-form-preserving action of ${\mathbb{Z}}/2{\mathbb{Z}}$ on the
``double octic''?

\subsection{Final comments on Table~\ref{thetable}}

As noted already, the authors do not know of any geometric
constructions of one parameter families of Calabi-Yau threefolds
with $h^{2,1}=1$ whose integral monodromy representations,
underlying their weight three ${\mathbb Z}$-VHS,  correspond to
the representations we found with $(m,a,t) = (1,4,1)$ or
$(2,4,2)$. In fact, we do not even know of candidate mirror
families with $h^{1,1} = 1$ and the correct geometric invariants.
Nevertheless, these integral monodromy representations do ``come
from geometry'' through restriction from multiparameter, higher
rank variations of the same weight.

The task of finding the 4D polytope (in case II) and 5D polytope
with NEF partition (in case I) with appropriate subloci is not
simple, but knowledge of the corresponding integral structures helps
to point the way.

The $(m,a)=(1,4)$ integral structure is unique, so it suffices to
merely find a sub-local system of the appropriate hypergeometric
type and then to prove that the restricted monodromy is actually
integral. The structure of the generalized hypergeometric series
$_4F_3(\frac{1}{12},\frac{5}{12},\frac{7}{12},\frac{11}{12})$
suggests looking for a complete intersection of bidegree $[2,12]$
in ${\mathbb{WP}}^5_{1,1,1,1,4,6}$. Unfortunately, there is not a
well-defined mirror of such a complete intersection (difficulties
with the singular locus). This is reflected in the fact that the
the Newton polytope for the weighted projective space is not
reflexive. We can correct this by considering instead the
reflexive polytope (provided by Kreuzer and Scheidegger) with
vertices given by $\{[e_1, e_3], [e_2, e_4, e_5,
-e_1-e_2-e_3-4e_4-6e_5, -2e_4-3e_5]\}$, where the brackets denote
the two pieces of the NEF partition. The complete intersection
Calabi-Yau threefold of bidegree $[2,12]$ in this Gorenstein Fano
toric variety has $h^{1,1} = 3$, but $h^{1,1}_{toric} = 2$. In
terms of the natural (torically defined) complex structure
coordinates $z_1, z_2$ for the polynomial deformations, the
restriction locus is just the locus $z_1 = 0$. Following Batyrev,
the generalized hypergeometric series for the holomorphic solution
in these coordinates is
$$\sum_{k,m} \frac{(2m)! (6k + 12m)!}{(3k + 6m)! (m!)^4 k! (2k + 4m)!} z_1^k
z_2^m \ ,$$ which upon restriction yields the desired series.  See
\cite{KKRS} for a discussion of the singular geometry and physics
of this subfamily.

The 4D polytope for the $(m,a)=(2,4)$ example is found by a similar
method.  Here one knows that the desired ${\mathbb{Z}}$-VHS has
$t=2$, so one expects if it arises by restriction that it will do so
from a family of Calabi-Yau threefolds with fundamental group
${\mathbb{Z}}/2{\mathbb{Z}}$. Also, the form of the hypergeometric
series $_4F_3(\frac{1}{8},\frac{3}{8},\frac{5}{8},\frac{7}{8})$
suggests looking for a hypersurface of degree $8$ in
${\mathbb{WP}}^4_{1,1,1,1,4}$.  Instead of taking this weighted
projective space as the starting point, since Kreuzer-Skarke have
completely classified the 473,800,776 equivalence classes of 4D
reflexive polytopes, one can imagine searching through these for
those hypersurfaces with nontrivial fundamental group.  In fact the
authors have done this, and a unique polytope was found with the
property that the hypersurface has fundamental group
${\mathbb{Z}}/2{\mathbb{Z}}$ and there is a sublocus of moduli on
which the hypergeometric series restricts to the desired form.

\end{document}